\documentclass[12pt]{amsart}
\usepackage{amssymb}

\def\C{\mathbb C}

\def\N{\mathbb N}

\begin{document}
\title{Escaping points of entire functions of small growth}
\author{P.J. Rippon and G.M. Stallard}
\date{}

\maketitle
\thanks{{\it Emails}: p.j.rippon@open.ac.uk,\; g.m.stallard@open.ac.uk}

{\it Abstract.} Let $f$ be a transcendental entire function and
let $I(f)$ denote the set of points that escape to infinity under
iteration. We give conditions which ensure that, for certain
functions, $I(f)$ is connected. In particular, we show that
$I(f)$ is connected if $f$ has order zero and sufficiently small
growth or has order less than 1/2 and regular growth. This shows
that, for these functions, Eremenko's conjecture that $I(f)$ has
no bounded components is true. We also give a new criterion
related to $I(f)$ which is sufficient to ensure that $f$ has no
unbounded Fatou components.

\section{Introduction}
\setcounter{equation}{0} Let $f$ be a transcendental entire
function and denote by $f^{n}, n \in \N$, the $n$th iterate of
$f$. The {\it Fatou set}, $F(f)$, is defined to be the set of
points, $z \in \C$, such that $(f^{n})_{n \in \N}$ forms a normal
family in some neighbourhood of $z$. The complement, $J(f)$, of
$F(f)$ is called the {\it Julia set} of $f$. An introduction to
the basic properties of these sets can be found in, for example,
[\ref{Berg}].

This paper concerns the {\it escaping set}
\[
 I(f) = \{z: f^n(z) \to \infty \mbox{ as } n \to \infty \},
\]
which was first studied for a general transcendental entire
function $f$ by Eremenko [\ref{E}]. He proved that

\begin{equation}
I(f) \neq \emptyset,
\end{equation}

\begin{equation}
J(f) = \partial I(f),
\end{equation}

\begin{equation}
I(f) \cap J(f) \neq \emptyset,
\end{equation}

\begin{equation}
\overline{I(f)} \mbox{ has no bounded components.}
\end{equation}
In particular, note that $I(f)$ is unbounded.

Eremenko conjectured that it may be possible to replace
$\overline{I(f)}$ with $I(f)$ in (1.4). This problem still
remains open although it has been shown to be true for certain
classes of functions -- see, for example, [\ref{Bar}] and
[\ref{RRRS}]. However, it was shown in~[\ref{RRRS}] that
points in $I(f)$ cannot necessarily be connected to infinity by a
curve in $I(f)$, answering another question of Eremenko.

For many of the functions considered in [\ref{Bar}] and [\ref{RRRS}], $I(f)$ consists of an infinite family of unbounded curves. On the other hand, Eremenko's conjecture is also true when $I(f)$ is connected, since $I(f)$ is always unbounded. In fact, $I(f)$ can be connected in surprisingly simple situations. For example, consider the function defined by $f(z)=z+1+e^{-z}$. This function was studied by Fatou [\ref{F}, Example 1] who showed that $F(f)$ consists of a single completely invariant domain, $U$ say, which is in $I(f)$. (Such a set $U$ is now known as a {\em Baker domain}.) For this function, $I(f)$ is connected because $U \subset I(f) \subset \overline{U}=\C$.

 In [\ref{RS}, Theorem 2] we proved the following result which shows that $I(f)$ is connected, and hence Eremenko's conjecture is true, in the more complicated situation where $F(f)$ has a multiply connected Fatou component. In the same paper, [\ref{RS}, Theorem~1], we showed that $I(f)$ always has at
least one unbounded component.

\newtheorem{1}{Theorem}
\begin{1}
Let $f$ be a transcendental entire function and suppose that
$F(f)$ has a multiply connected component $U$. Then $\overline{U}
\subset I(f)$ and $I(f)$ is connected.
\end{1}

{\it Remark.}
It is known [\ref{B84}, Theorem~3.1] that a multiply
connected component $U$ of $F(f)$ is a bounded wandering domain
and contains a Jordan curve $\gamma$ such that $f^n(\gamma) \to \infty$
and $f^n(\gamma)$ surrounds $0$ for sufficiently large $n \in
\N$.

In this paper we show that the escaping set is connected, and
hence Eremenko's conjecture is true, for other classes of
transcendental entire functions. In Section~2 we prove the
following generalisation of Theorem~1. This is the main result of
the paper. Note that $\widetilde{U}$ denotes the union of $U$ and
its bounded complementary components.

\newtheorem{2}[1]{Theorem}
\begin{2}
Let $f$ be a transcendental entire function and suppose that
there exists a bounded domain $U$ with $\partial U \subset I(f)$
and $U \cap J(f) \neq \emptyset$. Then

{\emph (a)} $\alpha_n = \partial \widetilde{f^n(U)}$ is in $I(f)$,
$\alpha_n \to \infty$ and $\alpha_n$ surrounds $0$ for
sufficiently large $n$;

{\emph (b)} $I(f)$ is connected.
\end{2}

Since $J(f) = \partial I(f)$, the hypothesis of Theorem~2 is
equivalent to the hypothesis that $\partial U \subset I(f)$ and
$U \cap I(f)^c \neq \emptyset$. We can thus describe this
hypothesis informally as there exists a `hole in $I(f)$'. It
follows from Theorem~2 that Eremenko's
conjecture is true whenever there is a hole in $I(f)$.

 We will show in the final section of the paper that there are many functions which
 satisfy the hypotheses of Theorem~2 but do not have any multiply
 connected Fatou components (so Theorem~1 does not apply to them).

We proved Theorem~1 in [\ref{RS}] by considering the following subset of
$I(f)$, which was introduced by Bergweiler and Hinkkanen in
[\ref{BH}]:
\[
A(f) = \{z: \mbox{there exists } L \in \N \mbox{ such that }
            |f^{n+L}(z)| > M(R,f^{n}), \mbox{ for } n \in \N\}.
\]
Here, for $r>0$,
\[
M(r,f) = \max_{|z|=r} |f(z)|
\]
and $R$ can be taken to be any value such that $R > \min_{z \in
J(f)}|z|$. We showed in~[\ref{RS}] that for any transcendental entire function $f$ the set
$A(f)$ is equal to the set
\[
  B(f) = \{z: \mbox{ there exists } L \in \N \mbox{ such that
  } f^{n+L}(z) \notin \widetilde{f^n(D)}, \mbox{ for } n \in \N \},
\]
where $D$ can be taken to be any open disc meeting $J(f)$, and also that $B(f)$ is connected under the hypotheses of Theorem~1; see [\ref{RS}, Theorem~2]. It follows from the `blowing up property' of the Julia set (see Lemma~2.1) that $B(f)$ is completely invariant.

In this paper, we introduce subsets of $B(f)$ of the form
\[
  B_D(f) = \{ z: f^n(z)\notin \widetilde{f^n(D)}, \mbox{ for } n \in \N \},
\]
where $D$ is any open disc meeting $J(f)$. Note that for any transcendental entire function $f$ all the components of $B(f)$ and $B_D(f)$ are unbounded; see [\ref{RS}, Theorem~1 and its proof]. In Section~3 we prove results concerning the structure of such sets which lead to various sufficient conditions for $B_D(f)$, $B(f)$ and $I(f)$ to be connected. Our first theorem in that section is the following.

\newtheorem{3}[1]{Theorem}
\begin{3}
Let $f$ be a transcendental entire function, let $D$ be an open
disc meeting $J(f)$, and suppose that there exists a bounded domain
$U$ with $\partial U \subset B_D(f)$
and $U \cap J(f) \neq \emptyset$. Then

{\emph (a)} $\beta_n= f^n(\partial U) \to \infty$ and, for
sufficiently large $n$, $\beta_n$ is in $B_D(f)$ and $\beta_n$ 
surrounds $0$;

{\emph (b)} $B_D(f)$, $B(f)$ and $I(f)$  are all connected.
\end{3}

We then show that there are several conditions that are
equivalent to the hypotheses of Theorem~3. This enables us to
prove the following result.

\newtheorem{4}[1]{Theorem}
\begin{4}
Let $f$ be a transcendental entire function, let $D$ be an open
disc meeting $J(f)$, and suppose that there exist Jordan curves
$\gamma_n$ such that, for all $n \in \N$,
\[
   \gamma_n \mbox{ surrounds } f^n(D)
\]
and
\[
  f(\gamma_n) \mbox{ surrounds the bounded component of }  \gamma_{n+1}^c.
\]
 Then there exists a bounded domain $U$ with $\partial U \subset
B_D(f)$ and $U \cap J(f) \neq \emptyset$. Thus $B_D(f)$, $B(f)$
and $I(f)$ are all connected.
\end{4}
Note that the hypotheses of Theorem~3 (and hence Theorem~4)
imply that the hypotheses of Theorem~2 are satisfied; that is, there exists a hole in $I(f)$.

We are aware of two different classes of functions
that satisfy the hypotheses of Theorem~4. Firstly, by the remark
following Theorem~1, any entire function with a multiply connected
component of the Fatou set will satisfy these hypotheses.
Secondly, many functions of order less than 1/2 have been shown
to satisfy conditions which, as we explain in Section~4, are
stronger than the hypotheses of Theorem~4. This gives various
results, including the following.

\newtheorem{5}[1]{Corollary}
\begin{5}
Let $f$ be a transcendental entire function and suppose that
either

{\it (i)} there exist $\epsilon \in (0,1)$ and $R>0$ such that
\begin{equation}
   \log \log M(r,f) < \frac{(\log r)^{1/2}}{(\log \log
   r)^{\epsilon}}, \; \mbox{ for } r>R,
\end{equation}
 or

{\it (ii)} the order of $f$ is $\rho < 1/2$ and
\begin{equation}
      \frac{\log M(2r,f)}{\log M(r,f)} \to c\; \mbox{ as } r \to \infty,
\end{equation}
where $c$ is a finite constant that depends on $f$.

Then, for any open disc $D$ meeting $J(f)$, the hypotheses of Theorem 4 hold, so $B_D(f)$, $B(f)$ and $I(f)$ are all connected.
\end{5}

Recall that the order of a function $f$ is defined to be
\[
   \rho = \overline{\lim_{r \to \infty}}\, \frac{\log \log M(r,f)}{\log
   r}\,.
\]

The functions of order less than 1/2 in Corollary 5 were
originally studied in connection with a different question,
associated with Baker [\ref{B81}], namely, whether a function of
order at most 1/2, minimal type, can have any unbounded Fatou
components.  It was shown in [\ref{S}, Theorem B and Theorem~C]
that the functions in Corollary 5 have no unbounded Fatou
components. Note that (1.5), which is not even satisfied by all
functions of order 0, is the best published growth condition on
$M(r,f)$ guaranteeing no unbounded Fatou components which does
not require some additional regularity condition such as (1.6).
See~[\ref{H}] for a survey article describing many other results
on this problem. Recently, we have shown that weaker conditions
than (1.5) and (1.6) are sufficient to ensure that a function of
order less than 1/2 has no unbounded Fatou components
(see~[\ref{SG}]), and Corollary~5 also holds under these weaker
conditions.

Finally, we give a criterion related to the escaping set which is
sufficient to ensure that a transcendental entire function has no unbounded Fatou
components.

\newtheorem{7}[1]{Theorem}
\begin{7}
Let $f$ be a transcendental entire function, let $D$ be an open
disc meeting $J(f)$, and suppose that there exists a bounded domain
$U$ with $\partial U \subset B_D(f)$ and $U \cap J(f) \neq
\emptyset$. Then $F(f)$ has no unbounded components.
\end{7}

Together with Theorem~4, this gives the following result.

\newtheorem{8}[1]{Theorem}
\begin{8}
Let $f$ be a transcendental entire function, let $D$ be an open
disc meeting $J(f)$, and suppose that there exist Jordan curves
$\gamma_n$ such that, for all $n \in \N$,
\[
   \gamma_n \mbox{ surrounds } f^n(D)
\]
and
\[
  f(\gamma_n) \mbox{ surrounds the bounded component of }  \gamma_{n+1}^c.
\]
Then $F(f)$ has no unbounded components.
\end{8}

{\it Remark.} The similarities between Theorems~3 and~6, and between Theorems~4 and~7, arise indirectly from their common methods of proof. It is natural to ask whether there is any direct relationship between the condition that $I(f)$ (or $B(f)$) is connected and the absence of unbounded
components of $F(f)$. Note, however, that the Fatou example $f(z)=z+1+e^{-z}$ mentioned earlier has an unbounded Fatou component and $I(f)$ is connected.

\section{Proof of Theorem~2}
\setcounter{equation}{0}

We begin this section by stating some existing results that will
be used in the proof of Theorem~2. Firstly, we will use the
following well-known property of the Julia set; see, for example,
[\ref{Berg}, Section~2] for a proof. This is often called the
`blowing up property' of the Julia set.

\newtheorem{2.1}{Lemma}[section]
\begin{2.1}
Let $f$ be a transcendental entire function, let $K$ be a compact set with $K \cap E(f) = \emptyset$ and let
 $U$ be an open neighbourhood of $z \in J(f)$. Then there exists
 $N \in \N$ such that
 \[
    f^n(U) \supset K,\quad \mbox{for all } n \geq N.
 \]
\end{2.1}

Here
\[
 E(f) = \{z: O^-(z) \mbox{ is finite} \},
\]
where
\[
 O^-(z) = \{w:f^n(w) = z, \mbox{ for some } n \in \N \}.
\]
The {\it exceptional set}, $E(f)$, contains at most one point.

We also use the following result; see [\ref{BH}] and [\ref{RS}].

\newtheorem{2.2}[2.1]{Lemma}
\begin{2.2}
Let $f$ be a transcendental entire function. The set $B(f)$ has the following properties:

{\it (a)} $J(f) = \partial B(f)$;

{\it (b)} $B(f)$ has no bounded components.
\end{2.2}

\begin{proof}[Proof of Theorem~2]

 Let $f$ be a transcendental entire
function and suppose that there exists a bounded domain $U$ with
$\partial U \subset I(f)$ and $U \cap J(f) \neq \emptyset$. Then
\[
 \alpha_n = \partial \widetilde{f^n(U)} \subset \partial f^n(U) \subset f^n(\partial U).
\]
So, since $I(f)$ is completely invariant, it follows from Lemma
2.1 that
\begin{equation}
\alpha_n \subset I(f), \; \alpha_n \to \infty \mbox{ and }
\alpha_n \mbox{ surrounds 0 for large enough } n.
\end{equation}

This proves part (a). We now show that $I(f)$ is connected.
We begin by showing that $B(f)$ belongs to one component of
$I(f)$. Suppose that $B_1$ and $B_2$ are two components of
$B(f)$. These are unbounded by Lemma~2.2(b) and so, by (2.1),
there exists $N \in \N$ such that

\[
 B_i \cap \alpha_n \neq \emptyset, \; \mbox{ for } i = 1,2 \mbox{ and } n \geq N.
\]

Since $\alpha_n \subset I(f)$, for each $n \in \N$, it follows
that $B_1$ and $B_2$ belong to the same component of $I(f)$. Therefore, there exists a component $I_0$ of $I(f)$ such that
\begin{equation}
B(f) \subset I_0 \; \mbox{ and } \; \alpha_n \subset I_0, \;
\mbox{ for all } n \geq N.
\end{equation}

Now suppose that $z_0 \in J(f) \cap I(f)$. It follows from Lemma
2.2(a) and (2.2) that $z_0 \in \partial B(f) \subset
\overline{I}_0$. Thus $I_0 \cup \{z_0\}$ is connected. Since $z_0
\in I(f)$ and $I_0$ is a component of $I(f)$, it follows that $z_0
\in I_0$. Thus
\begin{equation}
 J(f) \cap I(f) \subset I_0.
\end{equation}

To complete the proof we show that $F(f) \cap I(f) \subset I_0$.
Suppose that $V$ is a Fatou component in $I(f)$. (Note that a
Fatou component is either a subset of $I(f)$ or does not meet $I(f)$, by [\ref{Zip}, Theorem~3], for example.) We claim that there
exists a point in $\partial V \cap I(f)$. If $f^m(V)$ lies in a
multiply connected component $U$ of $F(f)$ for some $m \in \N$, then
it follows from Theorem~1 that $ f^m(\partial V) \subset \overline{U}\subset I(f)$. Thus $\partial V\subset I(f)$, by the complete invariance of $I(f)$.

Suppose on the other hand that $f^m(V)$ lies in a simply connected
component of $F(f)$ for all $m \in \N$. If $\partial V \subset
I(f)$, then the claim is proved. Otherwise, there exist $z_1 \in
\partial V$, $R>0$ and a sequence $(n_k)$ such that $|f^{n_k}(z_1)| < R$.
By (2.1) there exists $N \in \N$ such that
\begin{equation}
\alpha_N \mbox{ surrounds } \{z:|z| = R\} \; \mbox{ and } \;
B(0,R) \cap J(f) \neq \emptyset.
\end{equation}
Here we use the notation $B(w,r)=\{z:|z-w|<r\}$, where $w\in\C$ and $r>0$.

Since $V \subset I(f)$, there exists $k\in \N$ such that $f^{n_k}(V) \cap \alpha_N \neq \emptyset$. Since
$f^{n_k}(V)$ lies in a simply connected component of $F(f)$, there
are no curves in $f^{n_k}(V)$ surrounding points in $J(f)$. So, by
(2.4),
\[
\partial f^{n_k}(V) \cap \alpha_N \neq \emptyset.
\]
Thus $\partial f^{n_k}(V) \cap I(f) \neq \emptyset$. Since
$\partial f^{n_k}(V) \subset f^{n_k}(\partial V)$, we have
$\partial V \cap I(f) \neq \emptyset$ by the complete invariance
of $I(f)$.

So, whenever $V$ is a Fatou component in $I(f)$, we have
$\partial V \cap I(f)~\neq~\emptyset$. Since $\partial V \subset
J(f)$, it follows from (2.3) that $\partial V \cap I_0 \neq
\emptyset$. Now $V \cup (\partial V \cap I_0)$ is a connected
subset of $I(f)$ and so $V \subset I_0$. Hence $F(f) \cap I(f)
\subset I_0$. Together with (2.3), this shows that $I(f)$ is
connected. This completes the proof of Theorem~2.
\end{proof}

\section{Properties of the set $B_D(f)$}
\setcounter{equation}{0}

Let $f$ be a transcendental entire function and let $D$ be an
open disc meeting $J(f)$. Recall that
\[
  B_D(f) = \{ z: f^n(z)\notin \widetilde{f^n(D)}, \mbox{ for } n \in \N \}.
\]

In what follows we often use the property that, if $G$ is a
bounded domain, then
\[
  f(\widetilde{G}) \subset \widetilde{f(G)},
\]
which holds because if $\gamma$ is any simple closed curve in
$G$, then the image under $f$ of the inside of $\gamma$ lies
inside $f(\gamma)$ and hence in $\widetilde{f(G)}$.

\newtheorem{3.1}{Lemma}[section]
\begin{3.1}
Let $f$ be a transcendental entire function and let $D$ be an
open disc meeting $J(f)$. Then there exists $M \in \N$ such that
$f^m(B_D(f)) \subset B_D(f)$ for all $m \geq M$.
\end{3.1}

\begin{proof} Since $D \cap J(f) \neq \emptyset$, it follows from
the blowing up property, Lemma~2.1, that there exists $M \in {\bf
N}$ such that $\widetilde{f^m(D)} \supset D$ if $m \geq M$. Thus
\begin{equation}
\widetilde{f^{n+m}(D)} \supset \widetilde{f^n(\widetilde{f^m(D)})}
\supset \widetilde{f^n(D)},
\end{equation}
for each $n \in \N$, $m \geq M$. If $z \in B_D(f)$ then it
follows from (3.1) that, for each $n \in \N$,
\[
f^n(f^m(z)) = f^{n+m}(z)  \in  \C  \setminus
\widetilde{f^{n+m}(D)} \subset  \C  \setminus \widetilde{f^n(D)},
\]
 and hence $f^m(z)
\in B_D(f)$, for each $m \geq
M$. This completes the proof.
\end{proof}

\begin{proof}[Proof of Theorem~3]

Suppose that there exists a bounded domain $U$ with $\partial U
\subset B_D(f)$ and $U \cap J(f) \neq \emptyset$, and put
\[
 \beta_n = f^n(\partial U),\quad\text{for }n\ge 0.
\]
By Lemma~3.1, there exists $M\in\N$ such that $\beta_m\subset B_D(f)$ for $m\ge M$. In particular,
\[
 \beta_{n+M}=f^{n+M}(\partial U)\subset\C \setminus \widetilde{f^n(D)},
\]
for each $n\in \N$. Since $\partial f^{n+M}(U) \subset \beta_{n+M}$, part~(a) of Theorem~3 now follows from the blowing up property of the Julia set, Lemma~2.1.

We now recall that the proof of Lemma~2.2(a) in [\ref{RS}, proof
of Theorem~1] shows that there are no bounded components of
$B_D(f)$. Together with part~(a), this is sufficient to show that
$B_D(f)$ is connected. Similarly, since $B_D(f) \subset B(f)$ and
$B(f)$ has no bounded components (see Lemma~2.2(b)), it follows
from part~(a) that $B(f)$ is connected. Finally, since $\partial U\subset B_D(f)
\subset I(f)$, it follows from Theorem~2 that $I(f)$ is connected.
\end{proof}

The next result states that, unlike the sets $B(f)$ and $I(f)$,
the set $B_D(f)$ is always closed.

\newtheorem{3.2}[3.1]{Lemma}
\begin{3.2}
Let $f$ be a transcendental entire function and let $D$ be an
open disc meeting $J(f)$. Then $B_D(f)$ is closed.
\end{3.2}
\begin{proof} Since $D$ is open, we have $\widetilde{f^n(D)}$ is
open and hence $f^{-n}(\widetilde{f^n(D)})$ is open, for each $n
\in \N$. Thus
\[
  B_D(f)^c = \bigcup_{n \in \N}f^{-n}(\widetilde{f^n(D)})
\]
is open and hence $B_D(f)$ is closed.
\end{proof}

We now show that there are several properties of $B_D(f)$ that
are equivalent to the hypothesis on $B_D(f)$ given in Theorem~3.
These properties will be useful when proving that a given function satisfies
this hypothesis.

\newtheorem{3.3}[3.1]{Lemma}
\begin{3.3}
Let $f$ be a transcendental entire function and let $D$ be an open
disc meeting $J(f)$. The following are equivalent:

{\emph (a)} $B_D(f)^c$ has a bounded component;

{\emph (b)} there is a bounded domain $U$ with $\partial U
\subset B_D(f)$ and $U \cap J(f) \neq \emptyset$;

{\emph (c)} $\bigcup_{n \in \N}V_n$ is bounded, where $V_n$ is
the component of $f^{-n}(\widetilde{f^n(D)})$ that contains $D$.
\end{3.3}

\begin{proof} We begin by showing that (a) implies (b). So
suppose that $B_D(f)^c$ has a bounded component $U$. We know from
Lemma~3.2 that $U$ is open. Also $\partial U \subset B_D(f)$. Thus
\begin{equation}
  \partial f^n(U) \subset f^n(\partial U) \subset \C
  \setminus \widetilde{f^n(D)},
\end{equation}
for each $n \in \N$. Since $U \subset B_D(f)^c$, there exists $N
\in \N$ such that $f^N(U) \cap \widetilde{f^N(D)} \neq
\emptyset$. It follows from (3.2) that $f^N(U) \supset
\widetilde{f^N(D)}$ and so $f^N(U) \cap J(f) \neq \emptyset$. Thus
$U \cap J(f) \neq \emptyset$ and so (b) is true.

We now show that (b) implies (c). First note that if $V_n$ is the
component of $f^{-n}(\widetilde{f^n(D)})$ that contains $D$,
then  $\bigcup_{n \in \N}V_n \subset B_D(f)^c$. Also
\begin{equation}
V_n \subset V_{n+1}, \mbox{ for each } n \in \N,
\end{equation}
 since
\[
f^{n+1}(V_n) = f(f^n(V_n)) \subset f(\widetilde{f^n(D)}) \subset
\widetilde{f^{n+1}(D)}.
\]
Thus $\bigcup_{n \in \N}V_n$ is connected. If (b) is true then it
follows from Theorem~3(a) that there are no unbounded components
of $B_D(f)^c$ and so $\bigcup_{n \in \N}V_n$ is
bounded; that is, (c) is true.

Finally, we show that (c) implies (a). Put $V = \bigcup_{n \in
{\bf N}}V_n$ and suppose that (c) is true. Clearly $V \subset
B_D(f)^c$. We claim that $\partial V \subset B_D(f)$ and hence
(a) is true. To show this, we use proof by contradiction. So
suppose that $z_0 \in
\partial V$ and that $f^N(z_0) \in \widetilde{f^N(D)}$ for some
$N \in \N$. Since $\widetilde{f^N(D)}$ is open, there exists an
open neighbourhood $W$ of $z_0$ such that $f^N(W) \subset
\widetilde{f^N(D)}$ and hence
\begin{equation}
f^n(W) \subset \widetilde{f^n(D)}, \; \mbox{ for all } n \geq N.
\end{equation}
It follows from (3.3) that there exists $N_1 \geq N$ such that $W
\cap
\partial V_{N_1} \neq \emptyset$. Thus, by (3.4),
\[
 f^{N_1}(\partial V_{N_1}) \cap \widetilde{f^{N_1}(D)} \neq
 \emptyset.
\]
This, however, contradicts the fact that $V_{N_1}$ is a component of
$f^{-N_1}(\widetilde{f^{N_1}(V)})$ and so (a) is indeed true. This 
completes the proof of Lemma~3.3.
\end{proof}

 We are now in a position to prove Theorem~4.

%{\it Proof of Theorem~4.}
\begin{proof}[Proof of Theorem~4]

Let $f$ be a transcendental entire function and
let $D$ be an open disc meeting $J(f)$. Suppose that there exist
Jordan curves $\gamma_n$ such that, for all $n \in \N$,
\begin{equation}
   \gamma_n \mbox{ surrounds } f^n(D)
\end{equation}
and
\begin{equation}
  f(\gamma_n) \mbox{ surrounds the bounded component of }  \gamma_{n+1}^c.
\end{equation}

Now let $V_n$ be the component of $f^{-n}(\widetilde{f^n(D)})$
that contains $D$. We will show that $\bigcup_{n \in \N}V_n$ is
bounded. First note that, by the blowing up property, Lemma~2.1,
there exists $N \in \N$ such that
\begin{equation}
 \widetilde{f^N(D)} \supset D.
\end{equation}
We claim that
\[
 V_n \cap \gamma_N = \emptyset, \; \mbox{ for all } n \in \N.
\]
We prove this by contradiction. So, suppose that there exists $m
\in \N$ such that $V_m \cap \gamma_N \neq \emptyset$. Since $V_m
\supset D$, this implies that there exists a path $\sigma \subset
V_m$ joining a point in $D$ to a point in $\gamma_N$. By (3.7),
we have
\begin{equation}
f(D) \subset \widetilde{f^{N+1}(D)}
\end{equation}
and, by (3.5) and (3.6), $f(\gamma_N)$ surrounds the bounded
component of $\gamma_{N+1}^c$ and hence surrounds $f^{N+1}(D)$. So it
follows from (3.8) and (3.6) that $f(\sigma)$ joins a point in
$f(D)$ to a point in $\gamma_{N+1}$. Continuing this process
inductively, we find that $f^n(\sigma)$ joins a point in $f^n(D)$
to a point in $\gamma_{N+n}$. In particular, $f^m(\sigma)$ joins
a point in $f^m(D)$ to a point in $\gamma_{N+m}$. Since $\sigma
\subset V_m$ and $f^m(V_m) \subset \widetilde{f^m(D)}$, this
implies that
\begin{equation}
 \widetilde{f^m(D)} \cap \gamma_{N+m} \neq \emptyset.
 \end{equation}
By (3.7) we have $\widetilde{f^m(D)} \subset
\widetilde{f^{N+m}(D)}$ and so it follows from (3.9) that
\[
\widetilde{f^{N+m}(D)} \cap \gamma_{N+m} \neq \emptyset.
\]
This, however, contradicts (3.5) and so we must have
\begin{equation}
V_n \cap \gamma_N = \emptyset, \; \mbox{ for all } n \in \N.
\end{equation}
Since $V_n \supset D$, for all $n \in \N$, it follows from (3.5),
(3.7) and (3.10) that $V = \bigcup_{n \in \N} V_n$ lies in the
bounded component of $\gamma_N^c$ and is therefore bounded. So,
by Lemma~3.3, there exists a bounded domain $U$ with $\partial U
\subset B_D(f)$ and $U \cap J(f) \neq \emptyset$. Thus $B_D(f)$, 
$B(f)$ and $I(f)$ are all connected, by Theorem~3.
\end{proof}

\section{Functions of small growth}
\setcounter{equation}{0}

As mentioned earlier, it is conjectured that if $f$ is a transcendental entire function
of order at most $1/2$, minimal type, then $F(f)$ has no unbounded components.
This question was first studied by Baker [\ref{B81}] and has
since been studied by many other authors; see~[\ref{H}]. Many of the papers in this area
use the following result [\ref{S}, Lemma~2.7], which is
a generalisation of a result by Baker [\ref{B81}, proof of
Theorem~2].

\newtheorem{4.1}{Lemma}[section]
\begin{4.1}
Let $f$ be a transcendental entire function and suppose that
there exist sequences $R_n, \rho_n \to \infty$ and $c(n) > 1$
such that
\begin{enumerate}
\item $R_{n+1} = M(R_n,f)$,
\item $R_n \leq \rho_n \leq R_n^{c(n)}$,
\item $m(\rho_n,f) > R_{n+1}^{c(n+1)}$.
\end{enumerate}
Then $F(f)$ contains no unbounded components.
\end{4.1}

Here, for $r>0$,
 \[
  m(r,f) = \min_{|z|=r}|f(z)|.
\]

We now show that if $f$ satisfies the hypotheses of Lemma~4.1,
then $f$ also satisfies the hypotheses of Theorem~4, so the sets $B_D(f)$,
$B(f)$ and $I(f)$ are all connected.

\newtheorem{4.2}[4.1]{Lemma}
\begin{4.2}
Let $f$ be a transcendental entire function and suppose that
there exist sequences $R_n, \rho_n \to \infty$ and $c(n) > 1$
such that
\begin{enumerate}
\item $R_{n+1} = M(R_n,f)$,
\item $R_n \leq \rho_n \leq R_n^{c(n)}$,
\item $m(\rho_n,f) > R_{n+1}^{c(n+1)}$.
\end{enumerate}
Then, for each open disc $D$ meeting $J(f)$, there exist Jordan
curves $\gamma_n$ such that, for all $n \in \N$,
\begin{equation}
   \gamma_n \mbox{ surrounds } f^n(D)
\end{equation}
and
\begin{equation}
  f(\gamma_n) \mbox{ surrounds the bounded component of }  \gamma_{n+1}^c.
\end{equation}
Hence $B_D(f)$, $B(f)$ and $I(f)$ are all connected.
\end{4.2}

\begin{proof} Take $R_1$ such that $f(D) \subset B(0,R_1)$.
Then, if $R_{n+1} = M(R_n,f)$, we have $\widetilde{f^n(D)} \subset
B(0,R_n)$. If we put $\gamma_n = \{z:|z| = \rho_n\}$, then the
properties of the sequences $R_n$ and $\rho_n$ imply that (4.1)
and (4.2) are both satisfied. Thus $B_D(f)$, $B(f)$ and $I(f)$ 
are all connected, by Theorem~4.
\end{proof}

\begin{proof}[Proof of Corollary 5]

In [\ref{S}, Lemma 3.3, Lemma 4.2 and Lemma 4.5] we proved that functions satisfying (1.5) and functions of order $\rho < 1/2$ satisfying (1.6) must both satisfy the hypotheses of Lemma~4.1 (and so have no unbounded Fatou components); these proofs use estimates for $m(r,f)$ due to Baker [\ref{B58}] and Cartwright~[\ref{C}]. By Lemma~4.2, these functions also satisfy the hypotheses of
Theorem~4. This proves Corollary~5.
\end{proof}

\section{Unbounded Fatou components and the escaping set}
\setcounter{equation}{0}

In the previous section we stated that Lemma~4.1 forms the basis
of many proofs to show that certain functions of small growth have
no unbounded Fatou components. As we have seen, the
hypotheses of Lemma~4.1 are stronger than those of
Theorem~4 and hence stronger than those of Theorem~3. We now prove that the hypotheses of
Theorem~3 are in fact sufficient to ensure that a function has no
unbounded Fatou components -- that is, we prove Theorem~6. We use
the following results, the first of which is proved in
[\ref{Zip}, Theorem~3(b)], for example.

\newtheorem{5.1}{Lemma}[section]
\begin{5.1}
Let $f$ be a transcendental entire function, let $V$ be a Fatou component of $f$, and let $z_1, z_2 \in V$. If the Fatou components of $f$ which contain the forward orbits $f^n(z_i)$, $n=0,1,\ldots$, $i=1,2$, are all simply connected, then there
exists $A>0$ such that
\[
  |f^n(z_2)| \leq A(|f^n(z_1)| + 1), \; \mbox{ for all } n \in \N.
\]
\end{5.1}

Next we need a result about the size of the image of a large disc under
an iterate of a transcendental entire function. A weaker result of this type was given in [\ref{RS}, Lemma~2.2].

\newtheorem{5.2}[5.1]{Lemma}
\begin{5.2}
Let $f$ be a transcendental entire function. Then there exists
$R(f)>0$ such that
\[
  \widetilde{f^n(B(0,R))} \supset B(0,2M^n(R/2,f)), \; \mbox{ for
  all } n \in \N, R>R(f).
\]
\end{5.2}
\begin{proof} We first recall some facts from Wiman-Valiron
theory. For a detailed account, see~[\ref{Ha}], for example.

Let $f(z)=\sum_{n=0}^{\infty} a_n z^n$ be a transcendental entire function. For $r>0$, let $N(r)$ denote the integer $n\ge 0$ such that $|a_n|r^n$ is maximal and let $z_r$ satisfy $|z_r|=r$ and $|f(z_r)|=M(r,f)$. The main result of Wiman-Valiron theory states that if $\tau>1/2$, then there is a set $E\subset (0,\infty)$ such that $\int_E 1/t\,dt<\infty$ and \[
f(z)\sim \left(\frac{z}{z_r}\right)^{N(r)} f(z_r),\quad \mbox{for }z\in B(z_r,r/N(r)^{\tau}),
\]
as $r\to\infty$, $r\notin E$. In [\ref{E}, proof of Theorem~1], this result was used to prove that, for $r$ large enough and $r\notin E$, we have
\begin{equation}
f(B(z_r,5r/N(r))\supset \{z:e^{-4}M(r,f)<|z|<e^4 M(r,f)\}.
\end{equation}
Now $N(r)\to\infty$ as $r\to\infty$. Thus, for $R$ sufficiently large, there exists $r\in (R/2,R)$ such that~(5.1) holds and $B(z_r,5r/N(r))\subset B(0,R)$. Therefore,
\[
  \widetilde{f(B(0,R))} \supset B(0,2M(R/2,f)),
\]
for all sufficiently large $R$. The result then follows by induction.
\end{proof}

We are now in a position to prove Theorem~6.

\begin{proof}[Proof of Theorem~6]
Let $U$ be a bounded domain with $\partial U \subset B_D(f)$ and
$U \cap J(f) \neq \emptyset$, where $D$ is an open disc meeting
$J(f)$.

Suppose that $V$ is an unbounded component of the Fatou set.
Then $f$ has no multiply connected Fatou components, by the remark following Theorem~1. By Theorem~3(a), there exists $M \in \N$ such that
\[
 V \cap \beta_n \neq \emptyset, \; \mbox{ for all } n \geq M,
\]
where $\beta_n = f^n(\partial U)$, for $n\ge 0$, and so $V \subset I(f)$. We now
take $R_0$ so large that
\begin{equation}
  \beta_M \subset B(0,R_0)
\end{equation}
and $2M(R_0,f)> R(f)$, where $R(f)$ is as defined in Lemma~5.2.
Next we take $m \in \N$ such that
\begin{equation}
  \widetilde{f^{m+M}(D)} \supset B(0,2M(R_0,f)).
\end{equation}
(This is possible by the blowing up property, Lemma~2.1.)
Finally, we take $z_1 \in V \cap \beta_M$ and $z_2 \in V \cap
\beta_{m+M}$. Then, by (5.2),
\[
  |f^n(z_1)| \leq M(R_0,f^n), \; \mbox{ for all } n \in \N.
\]
Also, by (5.3) and Lemma~5.2, for all $n \in \N$,
\begin{eqnarray*}
  \widetilde{f^{n+m+M}(D)} & \supset & \widetilde{f^n(B(0,2M(R_0,f))}
   \supset B(0,2M^n(M(R_0,f),f))\\
   & = & B(0,2M^{n+1}(R_0,f)) \supset B(0,2M(M(R_0,f^n),f)).
\end{eqnarray*}
Thus, since $z_2 \in \beta_{m+M}$ and $\beta_0=\partial U \subset B_D(f)$,
so $f^n(z_2) \notin \widetilde{f^{n+m+M}(D)}$, we have
\[
  |f^n(z_2)| \geq 2M(M(R_0,f^n),f), \; \mbox{ for all } n \in \N.
\]
Hence
\[
  \frac{|f^n(z_2)|}{|f^n(z_1)|} \geq
  \frac{2M(M(R_0,f^n),f)}{M(R_0,f^n)} \to \infty \; \mbox{ as } n \to \infty,
\]
which contradicts Lemma~5.1. Therefore there are no
unbounded Fatou components.
\end{proof}

\section{Examples}
\setcounter{equation}{0}
Theorems~2,~3 and~4, and Corollary~5 give criteria on a transcendental 
entire function for the sets $B_D(f)$, $B(f)$ and $I(f)$ to be connected, where
$D$ is any open disc meeting $J(f)$.
Here we point out that there are many functions which satisfy one of the 
hypotheses of Corollary~5 (and hence satisfy the hypotheses of Theorems~2,~3
and~4) but which do not have multiply connected Fatou 
components (so Theorem~1 does not apply to them).
Indeed, we show that there are such functions of all orders $\rho$, 
$0\le\rho<1/2$. These are all new examples of
functions for which Eremenko's conjecture is true.

{\bf Example~1}\quad Bergweiler and Eremenko showed in [7] that
there are transcendental entire functions of arbitrarily small
growth for which the Julia set is the whole plane. In particular,
there exists such a function $f$ whose growth satisfies (1.5), and
hence has order 0. Thus, by Corollary~5(i), the sets $B_D(f)$, $B(f)$ 
and $I(f)$ are all connected. Clearly $f$ has no multiply connected Fatou components.

{\bf Example~2}\quad Baker [5] and Boyd [9] independently showed
that there are transcendental entire functions of arbitrarily
small growth for which every component of the Fatou set is simply
connected and every point in the Fatou set tends to~0 under
iteration. As before, there exists such a function $f$ whose
growth satisfies (1.5), and hence has order 0. Thus, by
Corollary~5(i), the sets $B_D(f)$, $B(f)$ and $I(f)$ are all connected. Clearly 
$f$ has no multiply connected Fatou components.

{\bf Example~3}\quad For $0<\rho<1/2$, consider the infinite product
\[
f(z)=c\prod_{n=1}^{\infty}\left(1+\frac{z}{n^{1/\rho}}\right),
\]
where $c>0$. Then, for $r>0$,
\[
M(r)=M(r,f)=f(r)\quad \text{and}\quad \log \frac{M(r)}{c}=r\int_0^{\infty}\frac{n(t)}{t(r+t)}\,dt,
\]
where $n(r)$ is the number of zeros of~$f$ in $\{z:|z|\le r\}$;
see~[\ref{T}, page~271]. Since the work of Wiman~[\ref{W}], such
functions have been known to have very regular behaviour. For
example, we have $n(r)=[r^{\rho}]$ so we can deduce by contour
integration that
\[
\log \frac{M(r)}{c}\sim r^{\rho}\frac{\pi}{\sin \pi\rho}\;\text{ as }r\to\infty.
\]
In particular, $f$ has order $\rho$ and~(1.6) holds, so the sets $B_D(f)$, $B(f)$ and $I(f)$ are all connected by Corollary~5(ii).

Similarly,
\[
r\frac{M'(r)}{M(r)} = r\int_0^{\infty}\frac{n(t)}{(r+t)^2}\,dt \sim r^{\rho}\frac{\pi\rho}{\sin \pi\rho}\;
\text{ as }r\to\infty,
\]
so there exist $C>1$ and $r_0>0$ such that
\begin{equation}
\frac{rf'(r)}{f(r)}\ge C\frac{\log f(r)}{\log r}\,,\quad\text{for }r>r_0.
\end{equation}
Now, for $r_0$ large enough, the interval $[r_0,\infty)$ is
invariant under $f$. Thus $[r_0,\infty)\subset J(f)$ by~(6.1) and [\ref{Zip}, Theorem~4(a)]. Hence $f$ has no multiply connected
Fatou components, by the remark following Theorem~1. Note that if
$c>0$ is small enough, then $f$ has at least one attracting fixed
point so $F(f)\ne \emptyset$.

{\it Remarks.}\quad 1. Example~3 can be generalised somewhat. Let
\[
f(z)=c\prod_{n=1}^{\infty}\left(1+\frac{z}{r_n}\right),
\]
where $c>0$, $r_n>0$. If we assume that $n(r)\sim \lambda r^{\rho}$, where $\lambda>0$ and $0<\rho<1/2$, then it can be shown as above that~(1.6) and~(6.1) hold (see~[\ref{T}, page~271] for the calculation related to $\log M(r)$). Thus the sets $B_D(f)$, $B(f)$ and $I(f)$ are all connected, and $f$ has no multiply connected Fatou components.

An example of such a function in closed form is given by
\[
f(z)=\frac12\left(\cos z^{1/4} + \cosh z^{1/4}\right)=1+\frac{z}{4!}+\frac{z^2}{8!}+\cdots.
\]
In this case $\rho=1/4$, $c=1$ and $f(-r)=\cos(\frac1{\sqrt{2}}\, r^{1/4})\cosh(\frac1{\sqrt{2}}\, r^{1/4})$, for $r>0$, so the zeros of $f$ are at $-4\pi^4(n-\frac12)^4$, $n=1,2,\ldots$. A family of examples of this type, with orders $\rho=2^{-m},\;m=2,3,\ldots,$ is given in~[\ref{Mar}, page~279].

2. Baker~[\ref{B85}] constructed examples of entire functions of all possible orders $\rho\ge 0$ which have multiply connected Fatou components. His examples were constructed using infinite products, the zeros being chosen in a more irregular manner than those in Example~3.

3. In [\ref{MF}] we remarked that it was plausible
that a transcendental entire function with no multiply connected Fatou
components should have infinitely many unbounded components of
$B(f)$. The above examples show that this is not the case.

\section*{References}
\begin{enumerate}

\item\label{B58} I.N. Baker. Zusammensetzungen ganzer Funktionen.
{\it Math. Z., 69 (1958), 121-163.}

\item \label{B81} I.N. Baker. The iteration of polynomials and
transcendental entire functions. {\it J. Aust. Math. Soc., Ser. A,
 30 (1981), 483--495.}

\item \label{B84} I.N. Baker. Wandering domains in the iteration
of entire functions. {\it Proc. Lond. Math. Soc., III. Ser., (3) 49 (1984),
563--576.}

\item \label{B85} I.N. Baker, Some entire functions with
multiply-connected wandering domains, {\it Ergodic Theory Dyn. Syst., 
5 (1985), 163--169.}

\item \label{B01} I.N. Baker. Dynamics of slowly growing entire
functions. {\it Bull. Aust. Math. Soc., 63 (2001), 367--377.}

\item \label{Bar} K.\ Bara\'nski.
Trees and hairs for some hyperbolic entire maps of finite order.
{\it Math.\ Z., 257 (2007), 33--59.}

\item \label{Berg} W. Bergweiler. Iteration of meromorphic functions.
{\it Bull. Am. Math. Soc., New Ser., 29 (1993), 151--188.}

\item \label{BE} W. Bergweiler and A.E. Eremenko. Entire functions
of slow growth whose Julia set coincides with the plane. {\it
Ergodic Theory Dyn. Syst., 20 (2000), 1577--1582.}

\item \label{BH} W. Bergweiler and A. Hinkkanen. On
semiconjugation of entire functions. {\it Math. Proc. Camb. Philos. Soc.,
126 (1999), 565--574.}

\item \label{Bo} D.A. Boyd. An entire function with slow growth
and simple dynamics. {\it Ergodic Theory Dyn. Syst., 22 (2002),
317--322.}

\item \label{C} M.L. Cartwright. {\it Integral functions.}
Cambridge Tracts in Mathematics and Mathematical Physics, No. 44,
Cambridge University Press, 1962.

\item \label{E} A.E. Eremenko. On the iteration of entire
functions. {\it Dynamical systems and ergodic theory,} Banach Cent. Publ.
 23 (Polish Scientific Publishers, Warsaw, 1989)
339--345.

\item \label{F} P. Fatou. Sur l'it\'{e}ration des fonctions 
transcendantes enti\`{e}res. {\it Acta Math. 47 (1926), 337-360.}

\item \label{Ha} W.K. Hayman. The local growth of power series: a
survey of the Wiman-Valiron method. {\it Can. Math. Bull. 17
(1974), 317--358.}

\item \label{H} A. Hinkkanen. Entire functions with bounded Fatou
components. To appear in {\it Transcendental dynamics and complex
analysis,}  Cambridge University Press, 2008.

\item\label{Mar} A.I. Markushevich. {\it Theory of functions of a
complex variable, Volume II,} Prentice Hall, 1965.

\item \label{Zip} P.J. Rippon and G.M. Stallard. On sets where
iterates of a meromorphic function zip towards infinity. {\it
Bull. Lond. Math. Soc., 32 (2000), 528--536.}

\item \label{RS} P.J. Rippon and G.M. Stallard. On questions of
Fatou and Eremenko. {\it Proc. Am. Math. Soc., 133 (2005),
1119--1126.}

\item \label{MF} P.J. Rippon and G.M. Stallard. Escaping points of
meromorphic functions with a finite number of poles. {\it J. Anal. 
Math., 96 (2005), 225-245.}

\item \label{SG} P.J. Rippon and G.M. Stallard. Functions of small
growth with no unbounded Fatou components. {\it Preprint.}

\item \label{RRRS} G.\ Rottenfu{\ss}er, J. R\"uckert, L. Rempe and D.\
Schleicher. Dynamic rays of bounded-type entire functions. {\it
Preprint.}

\item \label{S} G.M. Stallard. The iteration of entire functions
of small growth. {\it Math. Proc. Camb. Philos. Soc., 114 (1993),
43--55.}

\item \label{T} E.C. Titchmarsh. {\it The theory of functions}, Oxford University Press, 1939.

\item\label{W} A. Wiman. {\" U}ber die angen{\" a}herte Darstellung von ganzen Funktionen. {\it Ark. Mat., Astr. o. Fysik, 1 (1903), 105--111.}

\end{enumerate}

Department of Mathematics,\\
The Open University,\\
Walton Hall,\\
Milton Keynes MK7 6AA,\\
UK

E-mail: p.j.rippon@open.ac.uk; g.m.stallard@open.ac.uk

\end{document}